\documentclass[12pt]{amsart}
\usepackage{amssymb, amstext, amscd, amsmath}
\makeatletter
\def\@cite#1#2{{\m@th\upshape\bfseries%
[{#1\if@tempswa{\m@th\upshape\mdseries, #2}\fi}]}}
\makeatother
\theoremstyle{plain}
\newtheorem{thm}{Theorem}[section]
\newtheorem{cor}[thm]{Corollary}

\newtheorem{lem}[thm]{Lemma}
\theoremstyle{definition}
\newtheorem{rem}[thm]{Remark}
\newtheorem{defn}[thm]{Definition}

\newtheorem{eg}[thm]{Example}
\newcommand{\Prf}{\noindent\textbf{Proof.\ }}
\newcommand{\bx}{\hfill$\blacksquare$\medbreak}

\newcommand{\BH}{{\B(\H)}}

\newcommand{\ca}{\mathrm{C}^*}

\newcommand{\cenv}{\mathrm{C}^*_{\text{env}}}

\newcommand{\bbA}{{\mathbb{A}}}
\newcommand{\bbC}{{\mathbb{C}}}

\newcommand{\bbN}{{\mathbb{N}}}

  \newcommand{\A}{{\mathcal{A}}}
  \newcommand{\B}{{\mathcal{B}}}
  \newcommand{\C}{{\mathcal{C}}}

  \newcommand{\G}{{\mathcal{G}}}
\renewcommand{\H}{{\mathcal{H}}}
  \newcommand{\I}{{\mathcal{I}}}
  \newcommand{\J}{{\mathcal{J}}}
  \newcommand{\K}{{\mathcal{K}}}

  \newcommand{\N}{{\mathcal{N}}}

\renewcommand{\S}{{\mathcal{S}}}
  \newcommand{\T}{{\mathcal{T}}}

\renewcommand{\phi}{\varphi}

\newcommand{\fA}{{\mathfrak{A}}}

\newcommand{\fB}{{\mathfrak{B}}}

\newcommand{\fI}{{\mathfrak{I}}}
\newcommand{\fJ}{{\mathfrak{J}}}

\newcommand{\fM}{{\mathfrak{M}}}

\newcommand{\Alg}{\operatorname{Alg}}

\newcommand{\dirlim}{\varinjlim}
\newcommand{\dist}{\operatorname{dist}}

\newcommand{\Lat}{\operatorname{Lat}}

\newcommand{\ran}{\operatorname{Ran}}
\newcommand{\rank}{\operatorname{rank}}

\newcommand{\sumoplus}{\operatornamewithlimits{\sum\oplus}}
\newcommand{\supr}{\operatorname{sup}}
\newcommand{\Spec}{\operatorname{Spec}}

\newcommand{\nrep}{\operatorname{nrep}}

\begin{document}

%%%%%%%%%%%%%%%%%%%%%%%%%%%%%%%%%%%%%%%%%%%%%%%%%%%%%%%%%%%%
\title[]{Compact operators and nest representations of Limit Algebras}
\thanks{}
\author[E. Katsoulis]{Elias~Katsoulis}
\thanks{Katsoulis' research partially supported by a grant from
ECU}
\address{Math.\ Dept.\\East Carolina University\\
Greenville, NC 27858\\USA}
\email{KatsoulisE@mail.ecu.edu}
\author[J. Peters]{Justin~R.~Peters}
\address{Math.\ Dept.\\Iowa State University\\Ames, IA
50011\\USA}
\email{peters@iastate.edu}
\begin{abstract}
In this paper we study the nest representations \break $\rho: \A \longrightarrow \Alg \N$
of a strongly maximal TAF
algebra $\A$,
whose ranges contain non-zero compact operators. We introduce a particular
class of such representations, the
essential nest representations, and we show that their kernels coincide with the
completely meet irreducible ideals.
From this we deduce that there exist enough contractive nest representations,
with non-zero compact operators in their range,
to separate the points in $\A$.
Using nest representation theory, we also give a coordinate-free description of
the fundamental groupoid
for strongly maximal TAF algebras.

For an arbitrary nest representation $\rho: \A \longrightarrow \Alg \N$,
we show that the presence of
non-zero compact operators in the range of $\rho$
implies that $\N$ is similar to a completely atomic nest. If, in addition,
$\rho (\A )$ is closed
then every compact operator in $\rho (\A )$ can be approximated by sums of
rank one operators $\rho (\A )$.
In the case of $\bbN$-ordered nest representations, we show that $\rho ( \A
)$ contains finite rank operators iff
$\ker \rho $ fails to be a prime ideal.

\end{abstract}

\date{}
\maketitle
%%%%%%%%%%%%%%%%%%%%%%%%%%%%%%%%%%%%%%%%%%%%%%%%%%%%%%%%%%%%

%%%%%%%%%%%%%%%%%%%%%%%%%%%%%%%%%%%%%%%%
%%%%%%%%%%%%%%%%%%%%%%%%%%%%%%%%%%%%%%%%

\section{Introduction}

One of the central themes in the representation theory of $\ca$-algebras is
the consideration
of representations that contain non-zero compact operators in their ranges.
For instance, it is a deep result of Glimm~\cite{Gl} that characterizes the
type I $\ca$-algebras as those $\ca$-algebras
whose irreducible representations always contain non-zero compact operators
in their ranges. Other classes of $\ca$-algebras relating
to this line of thought are the CCR, GCR and residually finite
$\ca$-algebras, to mention but a few.

Recently there has been an interest in the representation theory of
non-selfadjoint algebras and in particular
triangular subalgebras of AF $\ca$-algebras~\cite{DK, DKPet, DHHLS, HPP00,
OPe}. In the case of non-selfadjoint
operator algebras, the irreducible representations are not always sufficient
to capture the structure of the
algebra and therefore appropriate substitutes have been introduced. A
representation $\rho$ of an
operator algebra $\A$ is said to be a \textit{nest representation} if the
closed invariant subspaces of $\rho (\A )$
are linearly ordered, i.e., form a \textit{nest}. Lamoureux \cite{La93}
coined the
term \textit{n-primitive ideal} for the kernel of a nest representation and
showed \cite{La93,La97} that in various contexts of
non-selfadjoint operator algebras
the n-primitive ideals play a role analogous to the primitive ideals in
$\ca$-algebras.
For instance, one can equip the set of n-primitive ideals with the
hull-kernel topology, and for every (closed, two-sided) ideal
$\J$ in the algebra, $\J$ is the intersection of all n-primitive ideals
containing $\J$.

In this paper we study nest representations whose
ranges contain non-zero compact operators.
The algebras under consideration are the \textit{strongly maximal TAF
algebras}, i.e., direct limits of sums
of upper triangular matrices so that the corresponding embeddings are
*-extendable and regular
\cite{PPW, MuS, Po}.
Perhaps surprisingly, there is little in the literature concerning
representations of strongly
maximal TAF algebras with compact operators in the range.
One reason may be that the ranges of the \textit{faithful} representations
of such algebras fail to contain
any compact operators, provided that the enveloping $\ca$-algebra is simple.
This applies to all
motivating examples of the theory, such as the standard, refinement and
alternation limit algebras. The situation changes
dramatically however if one allows the representation to have non-zero
kernel.
To begin with, every multiplicative linear form on the algebra
generates a representation on an one-dimensional Hilbert space and so its
range consists of rank one operators.
(We call these representations \textit{trivial}, for the obvious reasons.)
More interestingly, the standard limit algebra
has sufficiently many finite dimensional representations to separate the
points. (See \cite{LaS} for a characterization of the
residually finite TUHF algebras.)
It is not obvious at all however that other algebras, such as the refinement
or alternation limit algebras,
admit non-trivial nest representations with non-zero compact operators in
the range. It is this question
that initially motivated
the research in this paper. It turns out that this question is
related to a different open problem in representation theory which
we now describe.

An ideal $\J$ of a Banach algebra $\A$ is said to be \textit{ meet
irreducible} iff it cannot
be expressed as the intersection of a finite collection of ideals of $\A$
properly containing $\J$. The meet irreducible ideals are related to the
n-primitive ideals and this relationship was exploited in \cite{La93}. In
\cite{DHHLS} it was shown that for
a strongly maximal TAF algebra every
meet irreducible ideal is n-primitive. The converse question,
whether every n-primitive ideal is meet irreducible was left open in
\cite{DHHLS} and was answered affirmatively
in \cite{DKPet}.

Another problem that was left open in \cite{DHHLS} was the characterization
of completely meet
irreducible ideals using representation theory. An ideal $\J$ of a Banach
algebra $\A$ is called \textit{completely
meet irreducible} if, for any collection of ideals $\J_a$, $a \in \bbA$,
containing $\J$,
the relation $\J = \bigcap_{a  \in \bbA} \J_a$ implies $\J = \J_a$,
for some $a \in \bbA$. It is easily seen that this definition is
equivalent to the existence of a unique successor $\J^{+}$ of $\J$ in the
ideal lattice of $\A$.
The completely meet irreducible ideals have played a fundamental role in the
classification theory
of strongly maximal TAF algebras \cite{DHu, DHK, DPP}. Motivated by this,
and their study of nest representations,
the five authors of \cite{DHHLS} attempted a characterization
of the completely meet irreducible ideals of a strongly maximal TAF algebra
in terms of its nest representations (Corollary 5.7). In
\cite[Example~5.8]{DHHLS} it was realized that the order type of the nest
alone is not sufficient for such a characterization and therefore the
problem was left open.
The solution of this problem is contained in Theorem \ref{main} and
Corollary \ref{maincor}.
It turns out
that the missed point in \cite{DHHLS} was the use of finite rank operators.
Here we identify a
special class of nest representations whose ranges contain compact
operators, the essential nest representations (Definition
~\ref{CMI}),
and we show that the completely meet irreducible ideals coincide with the
kernels of the essential nest representations.
A key element in the proof of Theorem~\ref{main} is our recent
characterization of the meet irreducible
ideals as the kernels of the nest representations~\cite{DKPet}.

Utilizing the connection between compact operators and completely meet
irreducible ideals, we subsequently show that \textit{every} strongly
maximal TAF
algebra admits sufficiently many contractive representations, with compact
operators
in the range,
to separate the points and actually
capture the norm (Theorem \ref{norm}). This answers affirmatively the
question that
stimulated this work.

Motivated by the significance of nest representations with compact operators
in the range, we
undertake a general study of such representations in the second section of
the paper.
A first result shows that if the range of
a nest representation $\rho : \A \rightarrow \Alg \N$
of a strongly maximal TAF algebra $\A$ contains compact operators, then
the nest $\N$ is similar to a completely atomic one.
Having established that the nests involved in our considerations are
completely atomic (up to a similarity),
we subsequently consider the simplest possible case. In Theorem
\ref{N-ordered} we show that
the range of an $\bbN$-ordered nest representations $\rho : \A \rightarrow
\Alg \N$ contains finite rank operators if
and only if
$\ker \rho $ fails to be a prime ideal. (It is worth mentioning that the
proof of this result works
for any complex algebra.)
The second section
concludes with a result which shows that if the range of a nest
representation
$\rho : \A \rightarrow \Alg \N$ is closed then the compact operators in
$\rho ( \A ) $ form a closed ideal
in $\Alg \N$. This allows us to make use of the theory in \cite{ErdP} and
strengthen several of the results of this section, provided
that $\rho (\A )$ is norm closed.
Finally, note that
Example \ref{fewcompacts} shows that it is possible that the only non-zero
compact operators
in the range of a nest representation are rank one operators. Therefore we
cannot expect stronger
approximation results by utilizing strong operator topology. (Compare with
Corollary~\ref{rank one density}.)

In the final section of the paper we use representation theory
to give a coordinate-free description of the fundamental groupoid for a
strongly maximal TAF algebra $\A =\dirlim
( \A_i , \phi_i )$, i.e., a
description that does not depend on the ascending sequence of finite
dimensional subalgebras $\A_i$, $i \in \bbN$.
Our approach is influenced from the definition of the projective space from
Algebraic Geometry.

\begin{rem}
The results of Section~\ref{spectra} have influenced subsequent work on the
classification problem for non-selfadjoint operator algebras. In \cite{KaK}, the first author
and D. Kribs classify the quiver and free semigroupoid algebras of Muhly \cite{Muh}
and Kribs and Power \cite{KrP} using a variant of the
projective dual (Definition \ref{projdual}).
\end{rem}

%%%%%%%%%%%%%%%%%%%%%%%%%%%%%%%%%%%%%%%%%%%%%%%%%%%%%%%%%%%%%%%%%%%%%
%%%%%%%%%%%%%%%%%%%%%%%%%%%%%%%%%%%%%%%%%%%%%%%%%%%%%%%%%%%%%%%%%%%%%%%%%%%%%%%%%

\section{Compact operators, representations and completely meet irreducible
ideals}\label{existence}

In this section we characterize the completely
meet irreducible ideals as kernels of nest representations of a special kind.
Subsequently  we show that for every
strongly maximal TAF algebra there exist sufficiently many nest
representations, containing finite rank operators
in their ranges, to separate the points.

The following was also observed in \cite{HPP00} and will be repeatedly used
in the sequel.

\begin{lem} \label{onedimatom}
Let $\rho : \A \rightarrow \Alg \N$
be a nest representation of a strongly maximal TAF algebra $\A= \dirlim
(\A_i , \phi_i )$ on a Hilbert space $\H$.
Assume that $\rho$ is a $*$-representation on the diagonal of $\A$. If $p$
is an atom for
$\rho ( \A \cap \A^* )^{\prime \prime}$ then $\dim p(\H ) = 1$.
\end{lem}

\Prf If $g \in p( \H )$, then $[ \rho (\A )g ] \in \N$ and satisfies
\[
 [ \rho (\A )g ] = [g] \oplus M_g ,
 \]
for some closed subspace $M_g \subseteq p ( \H )^{\perp}$. Clearly any two
such subspaces $M_g , M_h$ cannot be
comparable unless $g = \lambda h$, $\lambda \in \bbC$, and the conclusion
follows.
\bx

If $\N$ is a nest then $0_+$ denotes the minimal non-zero element of $\N$
with the understanding
that if such an element does not exist then $0_+ = 0$. Similarly, $I_-$
denotes the largest element of
$\N$ not equaling $I$.

\begin{lem}  \label{rank one}
Let $\A= \dirlim (\A_i , \phi_i )$ be a strongly maximal TAF algebra and let
$\rho : \A \rightarrow \Alg \N$
be a nest representation of $\A$ on a Hilbert space $\H$. Assume that
$R \in \Alg \N $ is a non-zero
rank one operator so that
\[
R\rho (\A ) =  \rho (\A ) R = \bbC R .
\]
Then there exist non-zero vectors $g , h \in \H $ so that $0 _{+} = [ g ] $,
$I_{-}^{\perp} = [ h ] $ and
$R = g \otimes h$.
\end{lem}

\Prf Let $R = g \otimes h$ for suitable vectors $g , h \in \H $. Now $\rho
(\A )R  = \bbC R $
implies that $ [ \rho (\A )g ] = [g]$, i.e., the subspace $[ g ]$ is
invariant by $ \rho (\A ) $.
Since $\rho : \A \rightarrow \Alg \N$ is a nest representation, we obtain
$[g] \in \N$.
Therefore, $0 _{+} = [ g ] $. A similar argument shows that $I_{-}^{\perp} =
[ h ] $ and the
conclusion follows.
\bx

Corollary 5.7 of \cite{DHHLS} shows that for a certain class of nest
representations whose kernels are completely meet
irreducible ideals, $0$ has an immediate sucessor in the nest, and $I$ has
an immediate predecessor. In the following
theorem we require merely that the nest representation be a star
representation on the diagonal.

\begin{thm} \label{main}
Let $\A= \dirlim (\A_i , \phi_i )$ be a strongly maximal TAF algebra and let
$\rho : \A \rightarrow \Alg \N$
be a nest representation of $\A$ on a
Hilbert space $\H$. If $\rho$ is a $*$-representation on the diagonal
of $\A$, then the following are equivalent:
\begin{itemize}
\item[(i)] $\ker \rho$ is a completely meet irreducible ideal.
\item[(ii)] There exist non-zero vectors $g , h \in \H $ so that $0 _{+} = [
g ] $, $I_{-}^{\perp} = [ h ] $ and
$R = g \otimes h \in \rho ( \A )$.
\end {itemize}
\end{thm}

\Prf
Assume that $\J = \ker \rho$ is a completely meet irreducible ideal. Lemma
3.1 in \cite{DPP} implies
now that there exists a closed
ideal $\J^{+}$ containing $\J$ so that $\dim \left( \J^+ / \J \right ) = 1$.
Since closed ideals of
limit algebras are inductive, there exists a matrix unit $e \in \A_n$, $n
\in \bbN$, so that
$e \in \J^+ \setminus \J $. Furthermore, all but one of the subordinates of
$e$ in $\A_{i}$, $i\geq n $, are
mapped by $\rho $ to zero, or otherwise $\dim \left( \J^+ / \J \right )
\not= 1$. Therefore
we produce decreasing sequences $ \{ p_i \}_{i\geq n}$, $ \{ q_i \}_{i\geq
n}$ of diagonal matrix units,
 $p_i , q_i \in \A_i$, $i\geq n$, so that $p_i e q_i$ is a subordinate of
$e$ and
\begin{equation} \label{pq}
\rho ( p_i e q_i)  =  \rho (e).
\end{equation}
Note that the sequences
$ \{ \rho (p_i )\}_{i\geq n}$, $ \{ \rho (q_i ) \}_{i\geq n}$ are decreasing
sequences of selfadjoint projections and therefore
converge strongly to projections $p$ and $q$ respectively.

We claim that $p$ and $q$ are one-dimensional projections commuting with
$\N$.

First note that (\ref{pq}) shows that $p \rho(e) q =\rho(e)$ and so
both $p$ and $q$ are non-zero projections.
Now for a fixed $i \in \bbN$,
\[
p_{i+k} \A_i p_{i+k} = p_{i+k}p_i \A_i p_i p_{i+k} =\bbC p_{i+k},
\]
for all $k \in \bbN$.
Hence, $p \rho(\A_i) p = \bbC p $ and so $p \rho(\A ) p = \bbC p $. In
particular, $p$ is an atom for
$\rho ( \A \cap \A^* )^{\prime \prime}$ and so Lemma \ref{onedimatom} shows
that
$p$ is a rank one projection commuting with $\N$. A  similar argument proves
the claim for $q$

Since $p \rho(e) q =\rho(e)$ the above claim shows that $\rho (e )$ is a
rank one operator and so there
exist vectors $g , h \in \H $ so that $\rho(e) = g \otimes h$. Furthermore
notice that given any matrix unit
$a \in \A$, we have $ae \in \J^+$ and so either $ae = \lambda e + f$,
$\lambda \in \bbC, \ f \in \J$, or $ae \in \J
$ (and $f=0$). In any case, $\rho (a) \rho (e) = \lambda \rho (e)$ for some
$\lambda \in \bbC$. Hence, $ \rho(\A )  \rho(e)
= \bbC \rho (e)$ and similarly $\rho(e) \rho(\A )
= \bbC \rho (e) $, i.e., $\rho(e)$ satisfies the condition of Lemma
\ref{rank one} and the conclusion follows.

We now show that (ii) implies (i). A moment's reflection shows that $\bbC R
$ is an ideal in $\Alg \N$ and therefore
in $\rho ( \A )$. Hence, $\hat{\J} = \rho^{-1} ( \bbC R )$ is a norm closed
ideal of $\A$ containing $\ker \rho$.
Since $\dim \left( \hat{\J} / \ker \rho \right) = 1$, there exists $a \in
\A$ so that
\begin{equation} \label{dimension}
\hat{\J}= \bbC a + \ker \rho .
\end{equation}

Let $\I$ be any ideal properly containing $\ker \rho$. Then
\[
\ker \rho \subseteq \I \cap \hat{\J } \subseteq \hat{\J },
\]
with one of the above two containments being an equality because of
(\ref{dimension}).
However, $\ker \rho$ is the kernel of a nest representation and so Theorem
2.4 in \cite{DKPet} implies
that $\ker \rho$ is a meet irreducible ideal. Hence, $\I \cap \hat{\J } =
\hat{\J }$
and so $a \in \I$.
This shows that any closed ideal containing $\ker \rho$ also contains $a$
and so the
intersection of all such ideals properly contains $\ker \rho$. Conclusion:
$\ker \rho$ is a
completely meet irreducible ideal.
\bx

The statement of Theorem~\ref{main} suggests the following definition.

\begin{defn} \label{CMI}
A nest representation $\rho : \A \rightarrow \Alg \N$ of an operator algebra
$\A$ is
said to be \textit{essential} if the following conditions are
satisfied:
\begin{itemize}
\item[(i)] $\dim 0_+ = \dim I_{-}^{\perp} =1$.
\item[(ii)] The rank one operator from $I_{-}^{\perp}$ into $0_+$ belongs to
$\rho ( \A )$.
\end{itemize}
\end{defn}

We have arrived
at the desired characterization of completely meet irreducible ideals.

\begin{cor} \label{maincor}
A closed ideal of a strongly maximal TAF algebra is completely meet
irreducible ideal if and only if
it is the kernel of an essential nest representation.
\end{cor}

\Prf Observe that the proof of
(ii) implies (i) in Theorem~\ref{main} is applicable to any
essential nest
representation. (The
requirement that $\rho$ is a $*$-homomorphism
on the diagonal is used in the proof of the other direction.)
Therefore, the kernels of the essential nest representations are
completely meet irreducible ideals.
Conversely, any meet irreducible ideal $\J$ of a strongly
maximal TAF algebra $\A$ is the kernel of a nest representation $\rho$,
which is
a $*$-representation on the diagonal \cite{DHHLS}. If $\J$ is completely
meet irreducible then
the above theorem shows that $\rho$ is an essential nest representation and the
conclusion follows.
\bx

In \cite[Example~5.8]{DHHLS} it is shown that the above result fails if one
does
not require the range of the essential nest representation to contain the designated rank
one operator.
On the other hand, it is easy to see that for the refinement algebras, the
order type of the nest alone
implies the existence of the rank one operator in Definition~\ref{CMI}.
It would be interesting to know which are precisely the algebras which
satisfy that property.

%It turns out that for specific classes of limit algebras, the order type of
%the nest suffices for the existence
%of rank one operators in the range.
%
%\begin{thm} \label{orderpresmain}
%Let $\A= \dirlim (\A_i , \phi_i )$ be a strongly maximal TAF algebra so
%that the
%embeddings $\phi_i$ are order preserving. If $\rho : \A \rightarrow \Alg
%\N$
%is a nest representation of $\A$ on a
%Hilbert space $\H$, then the following are equivalent:
%\begin{itemize}
%\item[(i)] $\ker \rho$ is a completely meet irreducible ideal.
%\item[(ii)] $ \dim 0 _{+} = \dim I_{-}^{\perp} = 1 $.
%\end {itemize}
%\end{thm}

The following examples show that the finite rank operators contained in
the range of an essential nest representation may not be sufficient to describe the
range of the
representation in any reasonable topology.
To give the examples, we require the theory of $\ca$-envelopes as developed
in
\cite{DK}.

Let $\fA = \dirlim (\fA_i, \phi_i)$ be the enveloping $\ca$-algebra for
a TAF algebra $\A = \dirlim (\A_i, \phi_i)$ and let $\J \subseteq \A$ be
a closed ideal; let $\J_i := \J \cap \A_i$.
For each $i\ge1$, $\S_i$ denotes the collection of all diagonal projections
$p$
which are semi-invariant for $\A_i$, are supported on a single
summand of $\fA_i$ and satisfy $(p \A_i p) \cap \J = \{ 0 \}$. We form
finite dimensional C*-algebras
\[
 \fB_i := \sumoplus_{p \in S_i} \B(\ran p)
\]
where $\B(\ran p)$ denotes the bounded operators on $\ran p$;
of course, $\B(\ran p)$ is isomorphic to $\fM_{\rank p}$.
Let $\sigma_i$ be the map from $\fA_i$ into $\fB_i$
given by $\sigma_i(a) = \sum^\oplus_{p \in S_i} pap|_{\ran p}$.
The map $\sigma_i|_{\A_i}$ factors as
$\rho_i q_i$ where $q_i$ is the quotient map of $\A_i$ onto
$\A_i/ \J_i$ and $\rho_i$ is a completely isometric homomorphism of
$\A_i/ \J_i$ into $\fB_i$. Notice that $\fB_i$ equals the $\ca$-algebra
generated by $\rho_i ( \A_i/ \J_i )$.

We then consider unital embeddings $\pi_i$ of $\fB_i$ into $\fB_{i+1}$
defined as follows.
For each $q \in \S_{i+1}$ we choose projections $p \in \S_i$ which maximally
embed into $q$ under the action of $\phi_i$. This way, we determine
multiplicity one
embeddings of $\B (\ran p)$ into $\B (\ran q)$. Taking into account all such
possible embeddings,
we determine the embedding $\pi_i$ of $\fB_i$ into $\fB_{i+1}$.

Finally we form the subsystem of the directed limit $\fB= \dirlim (\fB_i,
\pi_i)$
corresponding to all summands which are \textit{never} mapped into a
summand $\B(\ran p)$ where $p$ is a maximal element of some $\S_i$.
Evidently this system is directed upwards.
It is also hereditary in the sense that if every image of a summand
lies in one of the selected blocks, then it clearly does not map into
a maximal summand and thus already lies in our system.
By \cite[Theorem~III.4.2]{D}, this system
determines an ideal $\fI$ of $\fB$.
The quotient $\fB' = \fB/ \fI$ is the AF algebra corresponding to the
remaining summands and the remaining embeddings; it can be expressed as a
direct
limit $\fB^{\prime} = \dirlim (\fB_{i}^{\prime}, \pi_{i}^{\prime})$, with
the understanding that
$\fB_{i}^{\prime}= \oplus_{j} \fB_{i\, j}$ for these remaining summands
$\fB_{i\, j}$
of $\fB_i$. It can be seen
that the quotient map is isometric on $A/ \J$ and that $\fB^\prime$
is the $\ca$-envelope of $A/ \J$.

\begin{eg} \label{fewcompacts}

There exists an infinite dimensional
essential nest representation
$\rho$ of a strongly maximal TAF algebra $\A$
so that the only compact operators in $\rho ( \A )$
are rank one operators.
{\medbreak}
Let $\A$ be the $2^{\infty}$-refinement
algebra and let $\J$ be the ideal of $\A$ determined by the sequence
\[
e^{(1)}_{1,2} , e^{(2)}_{1,3} , e^{(3)}_{1,5}, e^{(4)}_{1,9} , \dots
\]
of subordinates of the matrix unit $e^{(1)}_{1,2} \in \A_1$. Therefore,
the ideal $\A \cap \J_i$ is generated by all the matrix units in $\A_i$
except from
the ones in the "wedge"
 whose top right corner is $e^{(i)}_{1,2^{i-1}+1}$. Following the
construction of the
$\ca$-envelope described earlier, the $\ca$-envelope
$\fB^{\prime} = \dirlim (\fB_{i}^{\prime}, \pi_{i}^{\prime})$
of $\A / \J$ is given by the inductive limit
$$
M_{2}\bigoplus M_1\longrightarrow M_{3}\bigoplus M_{2}\longrightarrow
M_{5}\bigoplus M_{4}\longrightarrow \cdots
$$
where $M_{2^i}$ embeds in $M_{2^{i+1}}$ with the refinement embedding and
$M_{2^i +1}\oplus M_{2^{i}}$ embeds in $M_{2^{i+1}+1}$ with the embedding
that sends
\[
\left(
\begin{matrix}
 a_{1,1} & a_{1,2} & \ldots & a_{1,k+1}\cr a_{2,1} & a_{2,2} &
\ldots & a_{2,k+1}\cr \vdots & \vdots & \ddots & \vdots \cr a_{k+1,1} &
a_{k+1,2} & \cdots & a_{k+1,k+1}
\end{matrix}
\right)
\bigoplus
\left(
\begin{matrix}
 b_{1,1} & b_{1,2} & \ldots & b_{1,k}\cr b_{2,1}
& b_{2,2} & \ldots & b_{2,k}\cr \vdots & \vdots & \ddots & \vdots \cr
b_{k,1} &
b_{k,2} & \cdots & b_{k,k}
\end{matrix}
\right)
\]
to the matrix
\[
\left(
\begin{matrix}
 a_{1,1} & 0 & a_{1,2} & 0 & \ldots & a_{1,k} & 0 & a_{1,k+1}\cr
 0 & b_{1,1} & 0 & b_{1,2} & \ldots & 0 & b_{1,k} & 0 \cr
 a_{2,1} & 0 & a_{2,2} & 0 & \ldots & a_{2,k} & 0 & a_{2,k+1}\cr
 0 & b_{2,1} & 0 & b_{2,2} & \ldots & 0 & b_{2,k} & 0 \cr
 \vdots & \vdots & \vdots & \vdots & \ddots & \vdots & \vdots & \vdots \cr
 a_{k,1} & 0 & a_{k,2} & 0 & \ldots & a_{k,k} & 0 & a_{k,k+1}\cr
 0 & b_{k,1} & 0 & b_{k,2} & \ldots & 0 & b_{k,k} & 0 \cr
 a_{k+1,1} & 0 & a_{k+1,2} & 0 & \ldots & a_{k+1,k} & 0 & a_{k+1,k+1}
 \end{matrix}
\right) .
\]
We therefore have a commutative diagram
\[
\begin{CD}
\A_1 /J_1@> \phi_1>>  \A_2 / J_2 @> \phi_2>> \A_3 / J_3   @> \phi_3>>
\cdots
\\
@V\sigma_1VV    @V\sigma_2VV  @V\sigma_3VV  @.\\
M_{2}\oplus M_1 @> \pi_{1}^{\prime}>>
M_{3}\oplus M_{2}  @> \pi_{2}^{\prime}>> M_{5}\oplus M_{4}  @>
\pi_{3}^{\prime}>>
\cdots
\end{CD}
\]
where the vertical maps $\sigma_i$ are defined as follows. Given $a \in
\A_i$,
\begin{equation} \label{sigma}
\sigma_i (a+ \J)= \hat{a} \oplus \tilde{a},
\end{equation}
where $\hat{a}$ is the "wedge" of $a$ determined by $ a + \J_i$ and
$\tilde{a} $ results from $\hat{a}$ when its last row and column are
removed.

Consider the representation $\tau : \fB^{\prime} \rightarrow \BH$ of Theorem
2.4
in \cite{DKPet}. This is a faithful irreducible representation of
$\cenv (\A / \J )$  because $\J $ is meet irreducible. It also maps $\A /
\J$ densely
in a nest algebra $\Alg \N$. Let $\rho = \tau \circ \pi $, where $\pi: \A
\rightarrow  \A / \J$
is the quotient map.

Since $\J$ is a completely meet irreducible ideal, there exist
non-zero vectors $g , h \in \H $ so that $0 _{+} = [ h ] $, $I_{-}^{\perp} =
[ g ] $.
We claim that the only compact operators in $\rho ( \A ) = \tau (\A  / \J)$
are
of the form $g \otimes h^{\prime}$, $h^{\prime} \in \H$.

We identify first the rank one operators in $\tau \left(\bigcup _{i  \in
\Bbb N}
\fB_{i}^{\prime} \right)$. Let $e \in \bigcup _{i  \in \Bbb N}
\fB_{i}^{\prime}$. Since $\tau$ is an irreducible representation, $\tau(e)$
has rank one iff
\[
e \fB^{\prime} e = \bbC e ,
\]
i.e., $e$ has only one subordinate in any $\fB_{i}^{\prime}$. Therefore,
$e \in M_{2^{i+1}+1}$, for some $i \in \bbN$. If in addition $e$ belongs to
$\A / \J$ then (\ref{sigma}) shows that all columns of $e$ are $0$ except
the last
one. For such a $e$, $\tau (e)$ has the claimed form.

Finally, we look at $\A$ and consider $v \in \bigcup_{i \in \bbN} \A_i $ so
that
$\rho (v)$ is a rank one operator. Then $\pi(v)$ is mapped by $\tau$ to a
rank one
operator and therefore $\rho (v)= \tau ( \pi (v) )$ has the claimed form.
By Theorem \ref{mu compact}, sums of these operators approximate all
compacts in $\rho (\A )$
and the conclusion follows.
\end{eg}

With some additional work, one can obtain the following.

\begin{eg} \label{veryfewcompacts}
There exists an infinite dimensional
essential nest representation
$\rho$ of a strongly maximal TAF algebra $\A$
so that the only compact operators in
$\rho ( \A )$
are the scalar multiples of the rank one operator $R$ appearing in
Theorem~\ref{main}.
{\medbreak}
Let $\A$ be the $2^{\infty}$-refinement
algebra and this time let $\J$ be the ideal of $\A$ determined by the
sequence
\[
e^{(1)}_{1,2} , e^{(2)}_{2,4} , e^{(3)}_{3,7}, e^{(4)}_{6,14} , \dots
\]
of subordinates of the matrix unit $e^{(1)}_{1,2} \in \A_1$.
\end{eg}

We now show that for a strongly maximal TAF algebra $\A$, the nest
representations
whose ranges contain finite rank operators separate the points. First we
need a lemma.

\begin{lem} \label{cmi matrix unit}
Let $\A$ be a strongly maximal TAF algebra and let $e \in \A_i$, $i \in
\bbN$,
be a matrix unit. Then there exists a completely meet irreducible ideal $\J
\subseteq \A$ so that
$e \not \in \J$.
\end{lem}

\Prf Recall from \cite[Corollary 2.2]{DPP} that every element of the
spectrum of $\A$ can be identified
with a norm one linear functional which maps all matrix units onto $\{0, 1
\}$.
In \cite[Proposition 3.5]{DPP} it is shown
that for every meet irreducible ideal $\J$ there exists a unique functional
$\phi_{\J}$ in the spectrum
such that $\phi_{\J}$ annihilates $\J$ but does not annihilate $\J^{+}$.
The proof follows from that of \cite[Theorem 3.7]{DPP}. There it is shown
that if
$e \in \A_i$, $i \in \bbN$, is a matrix unit, then there exists a completely
meet
irreducible ideal $\J$ so that $\phi_{\J }(e) = 1$. Since $\phi_{\J}$
annihilates $\J$,
the conclusion follows.
\bx

The following result depends essentially on an application of Arveson's
distance formula.
 
\begin{lem}\cite{DK} \label{distance}
Let $\A= \dirlim ( \A_i , \phi_i)$ be a strongly maximal TAF algebra and let
$\J \subseteq \A$
be a closed ideal. If $a \in \A_i$, $i \in \bbN$, then
\[
\dist \left( a , \A_i \cap \J \right) = \dist \left( a , \J \right)
\]
\end{lem}

\Prf In \cite[Corollary 2.5]{DK} we show that the injection of $\A_{i} /
(\A_{i}\ \cap \J )$
into $\A_{i+1} / (\A_{i+1}\ \cap \J )$ is a complete isometry. Therefore, if
$a \in \A_i$, then
\[
\dist \left( a , \A_i \cap \J \right) = \dist \left( a ,  \A_{i+k} \cap \J
\right),
\]
for all $k \in \bbN$. The above equality passes to closed unions and since
$\J$
is inductive, the conclusion follows.
\bx

\begin{thm} \label{norm}
Let $\A= \dirlim ( \A_i , \phi_i)$ be a strongly maximal TAF algebra.
If $a \in \A$, then
\[
\| a \| = \supr \left\{\| \rho ( a )\|  \mid \rho \textit{ is a contractive
essential nest representation} \right\},
\]
\end{thm}

\Prf Let $\fA= \dirlim ( \fA_i , \phi_i)$ be
the enveloping $\ca$-algebra of $\A$. First we prove the result in the case
where $\fA$ is primitive.

Let $a_i \in \A_i$ so that $\| a - a_i \|\leq \epsilon$ and
so $\| a_i \|  \geq \| a \| - \epsilon$. Since $\fA$ is primitive there
exists a single summand
$\A_{j}^{(k_j )}$ in some $\A_j$, $j > i$, so that all summands of $\A_i$
are embedded in $\A_{j}^{( k_j )}$.
Let $e$ be the top right matrix unit in $\A_{j}^{(k_j )}$ and let $\J
\subseteq \A$ be a
completely meet irreducible ideal as in Lemma \ref{cmi matrix unit} so that
$e \not \in \J$.
By construction, none of the matrix units in $\A_i$ belongs to $\J$ or
otherwise $e \in \J$.
Therefore, $\A_i \cap \J = \{ 0 \}$ and so Lemma \ref{distance} shows that
\[
\dist ( a_i , \J) = \| a_i \| .
\]
Consider now an isometric nest representation $\hat{\rho}$ of $\A / \J $
\cite[Theorem 2.4]{DKPet}.
Then the induced representation $\rho: \A \rightarrow \Alg \N$ is a
contractive essential nest representation and
satisfies
\[
\| \rho (a_i)\|  = \dist ( a_i , \J)  = \| a_i \| \geq  \| a \| - \epsilon.
\]
Hence
\begin{eqnarray*}
\| \rho (a) \| & \geq & \|\rho (a_i )\| - \|\rho ( a - a_i )\|   \\
               & \geq & \| a \| -  2\epsilon.
\end{eqnarray*}
Since $\epsilon$ is arbitrary, the conclusion follows.

In order to prove the result in the general case, notice that there exists a
primitive ideal $\fJ$
of $\fA$ so that $\| a \| = \|a + \fJ \|$.
The natural map $\pi : \fA \rightarrow \fA / \fJ $ sends $\A$ on a strongly
maximal TAF algebra with primitive
enveloping $\ca$-algebra $\fA / \fJ $. From what we proved earlier it
follows that there exists a contractive essential nest representation $\hat{\rho}$ of
$\hat{\A }$ so that
\[
\| \hat{\rho} (a + \fJ ) \| \geq \|a + \fJ \| -\epsilon = \| a \| -\epsilon
\]
Considering the representation $\hat{\rho}\circ \pi$, the conclusion
follows.
\bx

%%%%%%%%%%%%%%%%%%%%%%%%%%%%%%%%%%%%%%%%%%%%%%%%%%%%%%%%%%%%%%%%%%%%%%%%%%%%
%%%%%%%%%%%%%%%%%%%%%%%%%%
%%%%%%%%%%%%%%%%%%%%%%%%%%%%%%%%%%%%%%%%%%%%%%%%%%%%%%%%%%%%%%%%%%%%%%%%%%%%
%%%%%%%%%%%%%%%%%%

\section{Structure for nest representations with compact operators in the
range}\label{existence}

In this section we study \textit{arbitrary} nest representations with
compact operators in the range.
The central result of the section is
the dichotomy of Theorem~\ref{N-ordered}; either the kernel of an
$\bbN$-ordered nest representation is
a prime ideal or else the range of the representation contains non-zero
finite rank operators. We also obtain some
approximation results, Theorem~\ref{mu compact} and Corollary~\ref{rank one
density}, which are similar in spirit
to those of \cite{Erd, KaM}.

We begin with a result that limits the types of nests involved in
representations
of strongly maximal TAF algebras with non-zero compact operators in the
ranges.
We need the following.

\begin{lem}
Let $\A= \dirlim (\A_i , \phi_i)$ be a strongly maximal TAF algebra
and $\rho : \A \rightarrow \Alg \N $ be a nest representation
on a Hilbert space $\H$. Assume that $\rho$ is a $*$-representation on the
diagonal $\A \cap \A ^* $
of $\A$.
If $\rho ( \A )$ contains non-zero compact operators then $\rho (\A \cap \A
^* )^{\prime \prime}$
has atoms.
\end{lem}

\Prf Notice that the elements of $\A$ that are mapped to compact operators
form a closed ideal in $\A$. Since closed ideals of limit algebras are
inductive, there exists an $n \in \bbN$
and a matrix unit $e \in \A_n$ so that $\rho ( e ) $ is a non-zero compact
operator. Now $e$ can be expressed
as a sum
\[
e =  \sum_{l} e^{(n+1)}_{l}
\]
of matrix units $e^{(n+1)}_{l} \in \A_{n+1}$ whose initial and final spaces
are orthogonal. Hence, the operators
 $\rho ( e^{(n+1)}_{l} )$ act on orthogonal subspaces of $\H$ and have
orthogonal ranges. Since their sum equals
 $\rho ( e ) $ there exists at least one of them, say $e_{n+1}$, so that $\|
\rho ( e_{n+1}) \| = \| \rho (e) \|$.
 Arguing as above we now produce a subordinate $e_{n+2} \in \A_{n+2}$ of
$e_{n+1}$ so that
 $\| \rho ( e_{n+2}) \| = \| \rho (e) \|$. Inductively, we define a sequence
$\{e_i \}_{i\geq n}$, $e_i \in \A_i$
 of matrix units each of which subordinates the previous one, so that $\|
\rho ( e_{i}) \| = \| \rho (e) \|$. for all $i$.
 Equivalently, we produce decreasing sequences $ \{ p_i \}_{i\geq n}$, $ \{
q_i \}_{i\geq n}$ of diagonal matrix units,
 $p_i , q_i \in \A_i$, $i\geq n$, so that $\| \rho ( p_i e q_i) \| = \| \rho
(e) \|$. Note that the sequences
$ \{ \rho (p_i )\}_{i\geq n}$, $ \{ \rho (q_i ) \}_{i\geq n}$ are decreasing
sequences of selfadjoint projections and therefore
converge strongly to projections $p$ and $q$ respectively. We will
show that the desired atom for
$\rho ( \A \cap \A^* )^{\prime \prime}$ is the projection $p$.

First we prove that $p$ is a
non-zero projection. Indeed, the sequence
$\{ \rho (p_i )\rho(e)\rho (q_i ) \}_{i\geq n}$ converges strongly to $p
\rho (e) q$. Since $\rho (e) $ is compact,
the sequence $\{ \rho (p_i )\rho(e)\rho (q_i  ) \}_{i\geq n}$ converges to $p \rho (e) q$
in norm. However,
\[
\| \rho (p_i )\rho(e)\rho (q_i ) \| =\| \rho (e_i ) \| =   \| \rho (e ) \|
\]
for all $i$ and so $\| p \rho (e) q \| =  \| \rho (e ) \| $, i.e., $p$ is a
non-zero projection.

Finally, notice that for a fixed $i \in \bbN$, $p_{i+k} \A_i p_{i+k} = \bbC
p_{i+k}$, for all $k \in \bbN$.
Hence, $p \rho(\A_i) p = \bbC p $ and so $p \rho(\A ) p = \bbC p $. In
particular, $p$ is an atom for
$\rho ( \A \cap \A^* )^{\prime \prime}$ and the conclusion follows.
\bx

\begin{thm} \label{completely atomic}
Let $\A$ be a strongly maximal TAF algebra
and $\rho : \A \rightarrow \Alg \N $ be a nest representation.
If $\rho ( \A )$ contains non-zero compact operators then $\N$ is similar
to a completely atomic nest. Furthermore, $\rho ( \A )$ is w*-dense
in $\Alg \N$.
\end{thm}

\Prf Recall that any bounded representation of $\A \cap \A ^* $ is
completely bounded~\cite[Theorem 8.7]{Paul},
and a completely bounded representation is similar to a completely
contractive one~\cite[Theorem 8.1]{Paul}.
Therefore $\rho$ is similar to a nest representation $\hat{\rho}$ so that
$\hat{\rho}$ is a $*$-representation on the diagonal of $\A$. Let $S$ be the
invertible operator
implementing that similarity, i.e., $\hat{\rho}(A) = S \rho (A) S^{-1}$, for
all $A \in \A$.
The previous lemma shows that
$\hat{\rho} (\A \cap \A ^* )^{\prime \prime}$ contains an atom. Proposition
3.5 in \cite{HPP00} shows now that
$\hat{\rho} (\A \cap \A ^* )^{\prime \prime}$ is a totally atomic von
Neumann algebra. By Lemma \ref{onedimatom}, the
atoms of $\hat{\rho} (\A \cap \A ^* )^{\prime \prime}$ are one dimensional
and so $\hat{\rho} (\A \cap \A ^* )^{\prime \prime}$
is a masa. Since $S\N S^{-1}$ commutes with $\hat{\rho} (\A \cap \A ^*
)^{\prime \prime}$, it is a completely atomic nest.

To prove the last sentence, notice that since $\hat{\rho} ( \A \cap \A^*
)^{\prime \prime}$  is a totally atomic masa in $\BH$, the
w*-closure
$\overline{\hat{\rho}  ( \A)}^{\, \textup{w}^* }$ of
$\hat{\rho}  ( \A)$ contains a masa and satisfies $\Lat \overline{\hat{\rho}
( \A)}^{\,\textup{w}^* } = S \N S^{-1}$. Therefore,
$\overline{\hat{\rho}  ( \A)}^{\, \textup{w}^* } = \Alg S \N S^{-1}$, by
\cite[Corollary 15.12]{Dav} and so
$\rho (\A)$ is w*-dense in $\Alg \N$.
\bx

Having established that the nests involved in our considerations are
completely atomic (up to a similarity),
we now undertake the simplest possible case. Recall that in \cite{DKPet} we
proved that
the kernel of any nest representation is a meet irreducible ideal. A
particular class
of meet irreducible ideals, the primitive ideals, were characterized in
\cite{DK} as the closed and prime ones.
The proof of the result below works for any complex algebra, not just
strongly maximal TAF algebras.
We therefore state it in that generality

\begin{thm} \label{N-ordered}
Let $\A$ be a complex algebra and let $\N$ be a maximal
nest which is ordered as $\bbN \cup \{ \infty \}$.
Let $\rho$ be a representation of $\A$ so that $\rho (\A)$ is w*-dense in
$\Alg \N$.
Then the following are equivalent:
\begin{itemize}
\item[(i)] $\rho ( \A ) $ contains non-zero finite rank operators.
\item[(ii)] $\ker \rho$ fails to be a prime ideal.
\end {itemize}
\end{thm}

\Prf Assume first that (ii) is valid and by way of contradiction assume that
 $\rho (\A)$ contains no non-zero finite rank operators.

Let $\J_1 , \J_2 $ be ideals in $\A$ so that $\J_1 , \J_2  \not\subseteq
\ker \rho$. Therefore the closed subspaces
$[ \rho (\J_i ) (\H ) ]$, $i = 1,2$, are non-zero invariant subspaces for
$\Alg \N$. Since $\rho ( \A )$
contains no finite rank operators, any non-zero operator in $\rho ( \J_i )$
has infinite dimensional range
and so we necessarily have
\[
[\rho (\J_1 )(\H ) ] = [ \rho (\J_2 ) (\H ) ] = \H.
\]
Let $a \in \J_1$ and $h \in \H$ so that $\rho(a)h \neq 0$. Since $[ \rho
(\J_2 ) (\H ) ] = \H$ there exist
sequences $\{ b_{n}^{(i)} \}_{n,i \in \bbN}$ and $\{ h_{n}^{(i)} \}_{n,i \in
\bbN}$ so that $\lim_{n \rightarrow \infty} \sum_{i}
\rho(b_{n}^{(i)})h_{n}^{(i)} = h$. Hence,
\[
\lim_{n \rightarrow \infty} \sum_{i} \rho(ab_{n}^{(i)})h_{n}^{(i)} = \rho (
a)h \neq0
\]
and therefore at least one of the terms $\rho(ab_{n}^{(i)})h_{n}^{(i)}$ is
non-zero. Hence $\J_1 \J_2  \not\subseteq \ker \rho$
and so $\ker \rho $ is prime, a contradiction.

We now prove that (i) implies (ii). Assume that $\rho (\A) $ contains a
non-zero finite rank operator $G$.
{\medbreak}
{\bf Claim 1.} $\rho( \A ) $ contains a non-zero finite rank operator $F$
with zero diagonal.
{\medbreak}
If $G$ has zero diagonal, then let $F = G$. Otherwise, the finite rank
operator $G$
has non-zero eigenvalues $\lambda$. The corresponding Riesz idempotents
$E(\lambda)$, $\lambda \in \bbC$
have finite
rank and belong to the algebra generated by $G$ and hence to $\rho ( \A)$.
Choose any such non-zero $E(\lambda)$
and let $S= \left[ E(\lambda)^* E(\lambda) -
(I - E(\lambda) )^* (I - E(\lambda))\right]^{1/2}$
be the normalizer of $E(\lambda)$. Thus
$ P \doteq S E(\lambda) S^{-1}$ is a selfadjoint projection that belongs to
the diagonal of the nest $S\N$.
(By a result of Larson $S\N$ is unitarily equivalent to $\N$ and therefore
an $\bbN$-ordered nest.)
Note that a similarity preserves the Jacobson Radical and therefore finite
rank operators with zero diagonal. Hence it is
enough to prove the claim for $S \rho ( \A ) S^{-1}$. We distinguish two
cases.

If $\dim P \geq 2$, then $PS \rho ( \A ) S^{-1}P$ is w*-dense and therefore
equal to
$P( \Alg S\N) P$ (Theorem~\ref{completely atomic}). Since $P (\Alg S\N ) P$
contains
rank one operators with zero diagonal, the same is true for
\[
PS \rho ( \A ) S^{-1}P = SE(\lambda)\rho ( \A )E(\lambda)  S^{-1} \subset S
\rho (\A ) S^{-1},
\]
which proves the claim in that case.

If $\dim P =1$, then $P \Alg S \N$ consists of rank one operators and has
\[
PS \rho ( \A ) S^{-1} = SE( \lambda)\rho ( \A ) S^{-1}
\]
as a dense subset. Hence there exists at least one non-zero rank one
operator in $SE( \lambda)\rho ( \A ) S^{-1}$
, say $R$. The desired finite rank operator will then be any operator of the
form
\[
F = R -t P
\]
for a suitable scalar $t$. The proof of the claim is complete.
{\medbreak}
Let $F$ be as in the Claim 1 and assume that $F$ has rank $n$. By
\cite[Lemma 1.2]{ErdP}, $F$ can be expressed as the sum
of $n$ rank one operators in $\Alg \N$ with zero diagonal. Hence there exist
projections $E_i \in \N$, $1 \leq i \leq n$,
and vectors $e_i \in E_{i}^{\perp}$, $f_i \in E_i$, $1 \leq i \leq n$, so
that
\[
F = \sum^{n}_{i=1} \, e_i \otimes f_i .
\]
Let $x$ be an element of $\A$ so that $\rho (x) = F$ and let $\I_x$ be the
smallest ideal of $\A$ containing
$x$.
{\medbreak}
{\bf Claim 2.} $( \I _x )^{n+1} \subseteq \ker \phi$.
{\medbreak}
It is enough to show that $( \I_F)^{n+1} = 0 $, where $\I_F$ denotes the
smallest ideal of $\Alg \N$ containing
$F$. Observe that $( \I_F)^{n+1}$ consists of finite sums of elements of the
form
\begin{equation} \label{product}
A_{i_1} \left(e_{i_1} \otimes f_{i_1} \right) A_{i_2} \left(e_{i_2} \otimes
f_{i_2} \right)   \dots
A_{i_{n+1}} \left(e_{i_{n+1}} \otimes f_{i_{n+1}} \right) A_{i_{n+2}},
\end{equation}
where $A_{i_j} \in \Alg \N$, for all $i_j$. At least two of the rank one
operators $e_{i_j} \otimes f_{i_j}$
coincide and therefore the term in (\ref{product}) reduces to
\[
A\left( e_{i_j} \otimes f_{i_j}\right) B \left( e_{i_j} \otimes
f_{i_j}\right) C
\]
with $A, B , C \in \Alg \N$. However,
\begin{eqnarray*}
A\left( e_{i_j} \otimes f_{i_j}\right) B\left( e_{i_j} \otimes
f_{i_j}\right)C & = &
A( E_{i_j}^{\perp}e_{i_j} \otimes f_{i_j}) B  \left( e_{i_j} \otimes
E_{i_j}f_{i_j}\right) C  \\
 & = &
A\left( e_{i_j} \otimes f_{i_j}\right)E_{i_j}^{\perp} B E_{i_j} \left(
e_{i_j} \otimes f_{i_j}\right) C  \\
               & = & 0
\end{eqnarray*}
since $B \in \Alg \N$. This proves the claim
{\medbreak}
Let $m$ be the least positive integer so that $( \I _x )^{m} \subseteq \ker
\phi$.
Since $\I_x \not \subseteq \ker \rho$, Claim 2 shows that
$m$ exists and $m \geq 2$. Consider the ideals $\I = \I_x$ and $\J =
(\I_x)^{m-1}$. Both are non-trivial ideals
not contained in $\ker \rho$
that satisfy $\I \J = ( \I _x )^{m} \subseteq \ker \phi$. Hence $\ker \rho$
is not a prime ideal, as desired.
\bx

All examples of $\bbN$-ordered nest representations of TUHF (triangular UHF)
algebras available in
the literature
fail to contain non-zero compact operators in the range. (The known
examples of
such representations with compact operators in the range occur for TAF
algebras
whose $\ca$-envelope
is not simple.)
The following example shows that such representations
do occur in a rather natural manner.

\begin{eg}
There exists a strongly maximal TUHF algebra $\A = \dirlim (\A_i , \phi_i )$
which
admits an $\bbN$-ordered nest representation with non-zero rank one
operators in the range.
{\medbreak}
Let $\A_i$, $i \geq 2$, be the collection of
all $3^i \times 3^i$-upper triangular matrices.
These are viewed as $3
\times 3$-block
upper triangular matrices of the form
\begin{equation} \label{matrix3}
\left(
\begin{matrix}
A_{1, 1} & A_{1, 2} & A_{1, 3} \\
0 & A_{2, 2} & A_{2, 3} \\
0 & 0 & A_{3, 3}
\end{matrix}
\right),
\end{equation}
where both the $(2, 2)$ and $(3, 3)$ entries
are $2^{i} \times 2^{i}$-upper triangular matrices. The embedding
$\phi_i$
maps (\ref{matrix3}) onto the matrix
\begin{equation} \label{matrix9}
\left(
\begin{matrix}
A_{1, 1}^{(3)} & A_{1, 2}^{(3)} & A_{1, 3}^{(3)} \\
0 & A_{2, 2}^{(3)} & A_{2, 3}^{(3)} \\
0 & 0 & A_{3, 3}^{(3)}
\end{matrix}
\right),
\end{equation}
where
\[ A_{k,l}^{(3)}= \left(
\begin{matrix}
A_{k,l} & 0 & 0 \\
0 & A_{k,l} & 0 \\
0 & 0 & A_{k,l}
\end{matrix}
\right) \]
is the threefold ampliation of $A_{k,l}$. (This is a variant of the familiar
block-standard embedding.)

Consider the collection $\J_i$ of all elements of $\A_i$
whose $(2, 2)$ and $(3, 3)$ entries in (\ref{matrix3}) are zero.
Evidently, $\J_i \subseteq \A_i$
is an ideal satisfying $\phi_i (\A_i ) \cap \J_{i+1} = \phi_i (\J_{i} )$.
Therefore the sequence $\{  \J_i \}_{i \in \bbN}$ determines a closed ideal
$\J$
so that $\J \cap \A_i = \J_i $, $i \in \bbN$. The quotient $\A / \J $ is the
strongly maximal
TAF algebra given by the inductive limit
$$
\begin{CD}
T_{2}\bigoplus T_2 @> \psi_{1} >> T_{4}\bigoplus T_{4} @> \psi_{2} >>
T_{8}\bigoplus T_{8} @> \psi_{8} >> \cdots
\end{CD}
$$
so that
$$
\psi_{i}
\left(
\begin{matrix}
A & 0 \\
0 & B
\end{matrix}
\right)
=
\left(
\begin{matrix}
\left(
\begin{matrix}
A & 0 \\
0 & B
\end{matrix}
\right)
& 0 \\
0 &
\left(
\begin{matrix}
B & 0 \\
0 & B
\end{matrix}
\right)
\end{matrix}
\right).
$$

Considering a variant of the Smith representation, one
sees that $\A / \J $ admits a $\bbN$-ordered nest representation
whose range contains all upper triangular compact operators and the
conclusion follows.
\end{eg}

The reader may have noticed that the proof of (i) implies (ii) in
Theorem~\ref{N-ordered} is valid for any
maximal nest $\N$. Therefore

\begin{cor}  \label{one direction}
Let $\A$ be a complex algebra, $\N$ be a maximal
nest and let $\rho$ be a representation of $\A$ so that $\rho (\A)$ is
w*-dense in $\Alg \N$.
If $\ker \rho$ is a prime ideal, then $\rho ( \A )$ does not contain
any non-zero finite rank operators.
\end{cor}

The refinement algebras are not semisimple and therefore the trivial ideal
$\{ 0 \}$ is not
prime in these algebras~\cite{DK}.
On the other hand, there exist many representations of the refinement
algebras whose
ranges do not contain any non-zero compact operators. This shows that the
converse
of Corollary~\ref{one direction} fails for nest representations
of order different than $\bbN$.

Note that in the proof of (i) implies (ii) in Theorem~\ref{N-ordered}, it
follows that if  $\rho ( \A )$ contains finite rank operators, then
$\A / \ker \rho $ cannot be semisimple and so $\ker \rho$ cannot be the
intersection of any collection of primitive ideals.

\begin{cor}
Let $\A$ be a strongly maximal TAF algebra and let $\rho$ be an
$\bbN$-ordered nest representation of $\A$ on a
Hilbert space $\H$. Then either $\ker \rho$ is a primitive ideal or $\ker
\rho$ is a meet irreducible ideal which cannot
be expressed as the intersection of primitive ideals. In that case, $\rho (
\A )$ contains non-zero finite rank operators.
\end{cor}

For the rest of this section, we specialize on nest representations with closed range.
Such representations exist in abundance. Indeed, Theorem 2.4 in \cite{DKPet} implies
that for any meet irreducible ideal $\J$ of a TAF algebra $\A$, there exists a nest representation
$\rho : \A \rightarrow \Alg \N$ with closed range so that $\ker \rho = \J$.
We begin with an approximation result that has already been used in
Example \ref{fewcompacts}

\begin{thm} \label{mu compact}
Let $\A= \dirlim ( \A_i , \phi_i)$ be a strongly maximal TAF algebra and let
$\rho : \A \rightarrow \Alg \N$ be a nest representation with closed range.
Let $e \in \bigcup_{ i \in \bbN}\A_i$ be an element of $\A$ so that
$\rho(e)$
is a compact operator. Then $\rho(e)$ can be expressed as
the sum of finitely many rank one operators in $\rho (\bigcup_{ i \in
\bbN}\A_i)$.
\end{thm}

\Prf Let $\J = \ker \rho$ and let $\hat{\rho} : \A / \J \rightarrow \Alg \N$
be the mapping
$\hat{\rho} (a + \J ) = \rho ( a)$, $a \in \A$.
Since $\rho ( \A )$ is closed, the Inverse Mapping Theorem implies that
$\hat {\rho}$ is an isomorphism of Banach spaces and so there exist
non-zero constants $c_1$ and $c_2$
so that
\begin{equation} \label{uniform}
c_1 \| a+ \J \| \leq \| \rho (a) \| \leq  c_2 \| a+ \J \| .
\end{equation}

It is enough to prove that for any matrix unit $e \in \A, \ \rho(e)$ can be
expressed as the sum of finitely many
rank one operators in $\rho (\bigcup_{ i \in \bbN}\A_i)$.
A matrix unit $e \in \A_n $ is said to be \textit{elementary} iff
$\rho(e) \not= 0$ and all but one of its subordinates in $\A_i$, $i \geq n$,
are annihilated by
$\rho$. An argument identical to that in the proof of Theorem~\ref{main}
shows that if $e$ is elementary
then $\rho (e)$ is a rank one operator.
 
We claim that $e \in \A_n$ is a finite sum of elementary operators.
Equivalently, we claim that
the number of subordinates of $e$ in $\A_i$, not annihilated by
$\rho$, is uniformly bounded. Clearly the claim proves the theorem.

By way of contradiction, assume that the number of subordinates of $e$ in
$\A_i$, not annihilated by
$\rho$, is not uniformly bounded. Then there exist two subordinates $e_1,
f_1 \in A_{i_1}$ of $e$ so that
$\rho (e_1), \rho(f_1) \not= 0$ and the number subordinates of $e_1$ in
$\A_i$, $i \geq i_1$,
 not annihilated by $\rho$ is not uniformly bounded. Consequently, there
exist two
 subordinates $e_2, f_2 \in A_{i_2}$ of $e_1$ so that
$\rho (e_2), \rho(f_2) \not= 0$ and the number subordinates of $e_2$ in
$\A_i$, $i \geq i_2$,
 not annihilated by $\rho$ is, once again, not uniformly bounded. This way
we produce a sequence
$\{ f_i \}_{i \in \bbN}$ of subordinates of $e$ so that $\rho ( f_i ) \not =
0 $. Furthermore, if
$q_i$ and $p_i$ are the initial and final spaces of $f_i$, $i \in \bbN$,
then
\[
p_i p_j = q_i q_j = 0
\]
for all $i, j \in \bbN$, with $i \not = \j$.

Since $\rho$ is a *-representation on the diagonal, $\{\rho ( p_i) \}_{ i
\in \bbN}$
and $\{ \rho (q_i ) \}_{ i \in \bbN}$ are sequences of mutually orthogonal
projections and therefore
converge strongly to $0$. Since $\rho (e)$ is a compact operator we obtain,
$$
\lim_{ i \in \Bbb N} \|\rho(  f_i) \|= \lim_{ i \in \Bbb N }\|\rho(p_i)
\rho(f) \rho(q_i) \|= 0.
$$
On the other hand,
\[
\dist(f_i , \A_i \cap \J) = \dist(f_i , \J)= 1
\]
and so (\ref{uniform}) implies that $\|\rho(f_i) \|\geq c_1$, which is
a contradiction and the conclusion follows.
\bx

\begin{cor} \label{rank one density}
Let $\A= \dirlim ( \A_i , \phi_i)$ be a strongly maximal TAF algebra and let
$\rho : \A \rightarrow \Alg \N$ be a nest representation with closed range.
Then every compact operator in $\rho (\A )$ can be approximated in norm by
sums of rank one
operators in $\rho ( \A )$. In particular, $\rho(\A )$ contains a non-zero
compact operator if and
only if it contains a non-zero rank one.
\end{cor}

\Prf The elements of $\A$ that are mapped to compact operators
form a closed ideal $\K$ in $\A$. Since the closed ideals of limit algebras
are inductive,
any such element $a \in \K$ can be approximated by elements in $\bigcup_{ i
\in \bbN}(\A_i \cap \K)$.
The previous theorem shows that $\rho \left(\bigcup_{ i \in \bbN}(\A_i \cap
\K) \right)$ consists of
sums of rank one operators in $\rho ( \A )$ and the conclusion follows.
\bx

\begin{thm} \label{ideal}
Let $\A= \dirlim ( \A_i , \phi_i)$ be a strongly maximal TAF algebra and let
$\rho : \A \rightarrow \Alg \N$ be a nest representation with closed range.
Then the compact operators in $\rho (\A)$ form a closed ideal of
$\Alg \N $.
\end{thm}

\Prf Corollary \ref{rank one density} shows that it is enough to check
products
between elements of $\Alg \N$ and rank one operators in $\rho (\A )$. Let
$e\otimes f$
be such a rank one operator and let $A \in  \Alg \N$. By
Theorem~\ref{completely atomic}
there exist nets $\{ B_i \}_i$ and  $\{ C^{*}_{j} \}_j$, $B_i , C_j \in
\rho( \A )$, converging to
$A$ and $A^*$ respectively in the strong operator topology. Hence, $\left\{
B_i\left (e\otimes f\right) \right\}_i$
and $\left\{ \left (e\otimes f\right) C_j\right\}_j$  are nets in $\rho
(\A)$ converging in norm to
$A (e\otimes f)$ and $(e\otimes f)A$ respectively. Since $\rho (A)$ is
closed, the conclusion follows.
\bx

Theorem~\ref{ideal} allows us to take advantage of the theory in
\cite{ErdP}. For instance,
if $\rho (\A)$ is closed then every rank $n$ operator in $\rho(\A)$ can be
written as the sum of $n$ rank one
operators of $\rho( \A)$ (\cite[Lemma 1.2]{ErdP}). This improves
Corollary~\ref{rank one density} in that
case. We can also improve on Theorems \ref{main} and \ref{N-ordered} as well.

\begin{cor}
Let $\A= \dirlim ( \A_i , \phi_i)$ be a strongly maximal TAF algebra and let
$\rho : \A \rightarrow \Alg \N$ be a nest representation with closed range.
If $\dim 0_{+}=\dim I_{-}^{\perp}=1$, then the following are equivalent:
\begin{itemize}
\item[(i)] $\rho(\A)$ contains non-zero compact operators.
\item[(ii)] $\ker\rho$ is a completely meet irreducible ideal of $\A$.
\end{itemize}
\end{cor}

\Prf Without loss of generality we
may assume that $\rho$ is a $*$-repre-
sentation
on the diagonal. If $\rho(\A)$ contains non-zero compact operators, then
Theorem~\ref{ideal} shows that the compact operators in $\rho (\A)$ form a closed ideal of
$\Alg \N $. Therefore, $\rho(\A)$ contains the rank one operator
$R = g \otimes h $ from item (ii) in Theorem \ref{main} and so
 $\ker\rho$ is a completely meet irreducible ideal of $\A$. Conversely,
if $\ker\rho$ is a completely meet irreducible ideal of $\A$, then
Theorem \ref{main} shows that $\rho(\A)$ contains a non-zero rank one
operator and the conclusion follows.
\bx

\begin{cor}
Let $\A= \dirlim ( \A_i , \phi_i)$ be a strongly maximal TAF algebra and let
$\rho : \A \rightarrow \Alg \N$ be a nest representation with closed range.
If $\N$ is a maximal
nest which is ordered as $\bbN \cup \{ \infty \}$,
then the following are equivalent:
\begin{itemize}
\item[(i)] $\rho(\A)$ contains non-zero compact operators.
\item[(ii)] $\ker\rho$ fails to be a prime ideal of $\A$.
\end{itemize}
\end{cor}

\Prf The result follows from Theorems \ref{rank one density}
and \ref{N-ordered}.
\bx

%%%%%%%%%%%%%%%%%%%%%%%%%%%%%%%%%%%%%%%%%%%%%%%%%%%%%%%%%%%%%%%%%%%%%%%%%%%%
%%%%%%%%%%%%%%%%%%%%%%%
%%%%%%%%%%%%%%%%%%%%%%%%%%%%%%%%%%%%%%%%%%%%%%%%%%%%%%%%%%%%%%%%%%%%%%%%%%%%
%%%%%%%%%%%%%%%%%%%%%%

\section{Spectra and representation theory for limit algebras} \label{spectra}

In this section we use representation theory to
describe an invariant for isometric isomorphisms of separable operator
algebras
that depends only on the algebra and its diagonal. Subsequently we show that
for
strongly maximal TAF algebras this invariant coincides with the fundamental
relation of Power~\cite{Po}.

\begin{defn}
If $\A$ is an operator algebra then $\nrep_{*} ( \A )$ denotes the collection of
all nest representations $\rho: \A \rightarrow \BH$ of $\A$ so that
$\rho$ is a $*$-representation on the diagonal $\A \cap \A^*$.
\end{defn}

Let $\rho \in \nrep_{*} (\A )$ and let $g , h \in \H$ be vectors so that the
subspaces $[g], [h]$ are atoms, i.e., minimal intervals,
 for $\Lat \rho (\A )$. (Note that if $\A $ is a strongly maximal TAF algebra,
 the existence of $[g], [h]$ implies that $\rho ( \A \cap \A^* )^{\prime \prime}$
 is a totally atomic masa and so $\N$ is totally atomic as well.)
 We define
 \[
 \omega^{(\rho)}_{g,h}(a) = \langle \rho (a) g , h \rangle, \quad a \in \A.
 \]
 The collection of all linear forms of the form $\omega^{(\rho)}_{g,h}$,
$\rho \in \nrep_{*} ( \A )$, $g , h \in \H$
is denoted as $\Omega_{ \A }$. We
do not wish to
distinguish between multiples of the same form. We therefore
define an equivalence relation $\sim$ on
$\Omega_{ \A }$ to mean
\[
\omega^{(\rho)}_{g,h} \sim \omega^{(\rho^{\prime})}_{g^{\prime},h^{\prime}}
\]
if and only if $\omega^{(\rho)}_{g,h} = \lambda
\omega^{(\rho^{\prime})}_{g^{\prime},h^{\prime}}$, for
some (non-zero) scalar $\lambda \in \bbC$. Consider now the quotient space
$\Omega_{ \A } / \! \sim$ and define a metric
$d$ on $\Omega_{ \A } / \! \sim$ by the formula
\[
d ( [ \omega_{1} ] , [ \omega_{2} ] )
= \sum_{n= 1}^{\infty} \frac{1}{2^n} \left| \,  \frac{\left| \, \omega_{1}(a_n) \, \right|}{\|\omega_{1}\|}  -
 \frac{\left|\, \omega_{2}(a_n) \, \right|}{\|\omega_{2}\|} \, \right|,
\]
where $\{a_n \}_{n=1}^{\infty}$ is a dense subset of $\A$. A moment's
reflection shows that
a sequence $\{[ \omega_{k} ] \}_{k=1}^{\infty}$ converges to some $[ \omega
] \in \Omega_{ \A } / \! \sim$ with respect to
the
metric $d$ iff
$$
\lim_{k \rightarrow \infty}  \frac{\left| \, \omega_{k}(a) \, \right|}{\|\omega_{k}\|}
= \frac{\left| \, \omega(a) \, \right|}{\|\omega\|}
$$
for all $a \in \A$. So even though the definition of the metric $d$
depends on the choice of the subset $\{a_n \}_{n=1}^{\infty}$, the topology
$\T$ determined
by $d$ is independent of that set.

\begin{defn}  \label{projdual}
If $\A$ is a separable operator algebra then the \textit{projective dual} of
$\A$ is the topological
space $\left( \Omega_{ \A } / \! \sim , \T \right)$.
\end{defn}

The adjective "projective" reflects the similarities between the above
construction and the
definition of the projective space in algebraic geometry.
Since isometric isomorphisms preserve the
diagonal, the projective dual is an invariant for isometric isomorphisms.

We compare now the projective dual of a strongly maximal TAF algebra to its
fundamental relation. Recall from
\cite{DPP} that the \textit{spectrum} or \textit{fundamental relation} of a
limit algebra $\A= \dirlim ( \A_i , \phi_i)$
is the set
\begin{equation} \label{Specdef}
\Spec (\A )= \left\{ \omega \in \A^{\#} \mid \| \omega \| =1 \mbox{ and }
\omega
\left( \bigcup_{i = 1 }^{\infty} \C_i \right) \subseteq \{ 0 , 1 \} \right\}
\end{equation}
equipped with the relative w*-topology as a subset of the dual $\A^{\#}$ of
the Banach space $\A$.
Here $\C_i$ denotes the natural matrix unit system associated with $\A_i$.
(The spectrum of $\A$ is also equipped with a partially defined operation
that makes it a
topological semigroupoid; we will come to these details later.)

\begin{lem} \label{omega}
If $\A= \dirlim ( \A_i , \phi_i)$ be a strongly maximal
TAF algebra, then $\Spec ( \A )\subseteq \Omega_{\A}$. Conversely,
if $\omega \in \Omega_{\A}$, then there exists a non-zero scalar $\lambda$
so that
$\lambda \omega \in \Spec (\A )$.
\end{lem}

\Prf Let $\omega \in \Spec ( \A )$. There exist decreasing sequences $\{ p_i
\}_{i = 1}^{\infty}$ and
$\{ q_i \}_{i = 1}^{\infty}$, $p_i , q_i \in \A_i$, $i \in \bbN$, of
diagonal projections and a matrix
unit $e \in \A_n $ so that
\[
\omega( q_i e p_i ) = 1
\]
and $\omega (f ) = 0$ for all other matrix units $f \in \A_i$, $i \geq n$.
The sequence
$\{ q_ie p_i \}_{i = 1}^{\infty}$ of subordinates of $e$ determines a path
$\Gamma$
on the Bratelli diagram for
the enveloping $\ca$-algebra $\fA= \dirlim ( \fA_i , \phi_i)$ of $\A$.
Let $\fJ$ be the ideal consisting of all
summands of $\fA= \dirlim ( \fA_i , \phi_i)$ that are never mapped into
$\Gamma$
and let $\pi: \fA \rightarrow \fA / \fJ$ be the quotient map. Then $\pi (
\A)$ is a strongly
maximal TAF algebra of a primitive $\ca$-algebra. Furthermore, the
decreasing sequence
$\{ \pi( p_i) \}_{i = 1}^{\infty}$ of diagonal projections determines a
multiplicative linear form $p$
on $\pi(\A) \cap \pi(\A )^*$, whose orbit under the action of $\fA / \fJ$ is
dense in the Gelfand
spectrum of $\pi(\A) \cap \pi(\A )^*$. We now apply the construction of
\cite[Proposition II.2.2.]{OPe}
so that the role of $[x_0]$ in that proof is played by our $p$. Therefore we
obtain a faithful
nest representation
$\hat{\rho}$ of $\pi (\A )$ so that $\hat{\rho}$ is a *-representation on
the
diagonal and $\Lat  \hat{\rho}\left(\pi(\A) \right)$ is a maximal, totally
atomic nest. In addition,
the construction of \cite[Proposition II.2.2.]{OPe} implies that the
sequences
the sequences $\{ \hat{\rho}\left(\pi (p_i) \right) \}_{i = 1}^{\infty}$ and
$\{ \hat{\rho}\left(\pi (q_i ) \right) \}_{i = 1}^{\infty}$ converge to one
dimensional projections
with ranges $[g]$ and $[h]$ respectively. If $\rho = \hat{\rho}\circ \pi$
then it is easy to see that
$\omega = \omega^{(\rho)}_{g,h}$ and the conclusion follows.

Conversely, let $\omega \in \Omega_{\A}$ and let $\rho \in \nrep_{*} ( \A )$, $g,
h \in \H$ so that
$\omega = \omega^{(\rho)}_{g,h}$. There exists a matrix unit $e \in \A_n$ so
that
$\omega(e) \not= 0$ and so let $\lambda$ be a scalar so that
$\omega(e) = \lambda^{-1}$.
Since $\rho$ is a *-representation on the diagonal and $[g],[h]$ are atoms
for
$\rho ( \A \cap \A^* )^{\prime \prime}$, $\omega$ annihilates all but one of
the matrix units of
$\A_i$, $i \geq n$; let us denote as $e_i$ that matrix unit in $\A_i$. A
moment's reflection shows that
$e_i$ is a subordinate of $e$ and so $\omega(e_i) = \lambda^{-1}$. It
is clear now that
$\lambda \omega$ satisfies the requirements of (\ref{Specdef}) and the proof
follows.
\bx

Note that Theorem~\ref{omega} establishes the existence of an 1-1 and onto
map $\Phi: (\Spec ( \A ), \mbox{w}^* )
\rightarrow ( \Omega_{ \A } / \! \sim , \T ) $. Furthermore, this map is
easily seen to be bicontinuous. Hence

\begin{cor} \label{bicont}
The topological spaces $(\Spec ( \A ), \mbox{w}^* )$ and $( \Omega_{ \A } / \! \sim , \T )$
are homeomorphic.
\end{cor}

In the case of a strongly maximal TAF algebra generated by its order
preserving normalizer, Theorem~\ref{main}
shows that the essential nest representations are sufficient to
describe
the fundamental relation. In the general case on obtains a dense subset of
$( \Omega_{ \A } / \! \sim , \T )$,
as \cite[Theorem 3.7]{DPP} shows.

The fundamental relation is also equipped with partially defined operation
that makes it a
topological semigroupoid. We now indicate how this operation materializes in our
representation.

Let $\A$ be a separable operator algebra and let
$\omega^{(\rho)}_{g,h} , \omega^{(\rho^{\prime})}_{g^{\prime},h^{\prime}}
\in \Omega_{\A}$. If $\rho = \rho^{\prime}$ and $[g] = [h^{\prime}]$ then
we define
\begin{equation} \label{composition}
\omega^{(\rho)}_{g,h} \circ \omega^{(\rho^{\prime})}_{g^{\prime},h^{\prime}}
=
\omega^{(\rho)}_{g^{\prime},h}.
\end{equation}
It is not clear at all that (\ref{composition}) establishes a well
defined operation on $\Omega_{\A}$. Nevertheless, this is the case for strongly maximal
TAF algebras, as the following result shows.

\begin{thm} \label{final}
Let $\A$ be a strongly maximal TAF algebra. Then the triple
$( \Omega_{ \A } /\!\sim , \T , \circ)$ forms a topological semigroupoid
which is isomorphic to the fundamental relation of $\A$.
\end{thm}

\Prf If $\omega^{(\rho)}_{g,h} \in \Phi( \omega )$, then the last paragraph in the proof of
Lemma~\ref{omega} shows that $\omega^{(\rho)}_{g,g} = s(\omega)$
and $\omega^{(\rho)}_{h,h} = r (\omega)$, where $s(\cdot)$ and $r (\cdot )$ are
the source and range maps~\cite{DPP}. Hence, if
$\omega^{(\rho)}_{g,h} \in \Phi (\omega)$ and
$\omega^{(\rho^{\prime})}_{g^{\prime},h^{\prime}} \in \Phi ( \omega^{\prime})$
are composable, then $r(\omega^{\prime}) = s(\omega)$ and
so $\omega$ and $\omega^{\prime}$ are composable.
In addition, $\omega^{(\rho)}_{g^{\prime},h} \in \Phi (\omega \circ \omega^{\prime})$. This
shows that $\circ$ is a well defined
operation on $\Omega_{ \A } / \!\sim $ and that $\Phi^{-1}$ respects that operation.

We now show that $\Phi$ also respects the composition on $\Spec (\A )$.
Let $\omega , \omega^{\prime} \in \Spec ( \A )$ with $r(\omega^{\prime}) =
s(\omega)$. Then  $r(\omega) , s(\omega), r(\omega^{\prime})$ and $s(\omega^{\prime})$
belong to the same orbit under the action of $\A$ on the Gelfand space of the diagonal.
Therefore the construction of \cite[Proposition II.2.2.]{OPe} produces a common
representation $\rho$ for both $\omega$ and $\omega^{\prime}$ so that
$\omega^{(\rho)}_{g,h} \in \Phi (\omega)$ and
$\omega^{(\rho^{\prime})}_{g^{\prime},h^{\prime}} \in \Phi ( \omega^{\prime})$.
Since $r(\omega^{\prime}) =
s(\omega)$, the construction of \cite[Proposition II.2.2.]{OPe} also guarantees
that $g = h^{\prime}$. Hence $\omega^{(\rho)}_{g,h}$ and
$\omega^{(\rho^{\prime})}_{g^{\prime},h^{\prime}}$ are composable. In addition,
\[
\Phi (\omega \circ \omega^{\prime})
=[\omega^{(\rho^{\prime})}_{g^{\prime},h}]
=[\omega^{(\rho)}_{g,h} ] \circ [\omega^{(\rho^{\prime})}_{g^{\prime},h^{\prime}}]
= \Phi (\omega) \circ \Phi(\omega^{\prime})
\]
and so $\Phi$ also respects the composition. From this it is easy to see that
the operation $\circ$ on $\Omega_{ \A } /\! \sim $ is associative. The bicontinouity of
$\Phi$ has been established earlier and the continuity of $\circ$ is easy to prove.
\bx

Notice that Theorem~\ref{final} gives yet another proof of the invariance of the fundamental
relation under isometric isomorphisms. The triple
$( \Omega_{ \A } /\!\sim , \T , \circ)$ can be calculated for other operator algebras
as well,
including various semicrossed products. It seems that in the case where
$\A$ is the algebra of an r-discrete, principal semigroupoid $\G$~\cite{MuS},
one should be able to relate the triple
$( \Omega_{ \A } /\!\sim , \T , \circ)$ to the semigroupoid $\G$. We plan to pursue
these directions elsewhere.
%%%%%%%%%%%%%%%%%%%%% References %%%%%%%%%%%%%%%%%%%%%


\begin{thebibliography}{99}

\bibitem{Dav} K. Davidson,
\textit{Nest Algebras},
Pitman Research Notes in Mathematics Series, \textbf{191}, 1988.

\bibitem{D} K. Davidson,
\textit{$\ca$-algebras by example},
Fields Institute Monographs, American Mathematical Society, 1996.

\bibitem{DK} K. Davidson and E. Katsoulis,
\textit{Primitive limit algebras and $\ca$-envelopes},
Adv.\ Math. \textbf{170}, 2002, 181--205.

\bibitem{DKPet} K. Davidson, E. Katsoulis and J. Peters,
\textit{Meet Irreducible Ideals and Representations of Limit Algebras},
J. Funct. Anal. \textbf{200}, 2003, 23--30.

\bibitem{DHHLS} A. Donsig, A. Hopenwasser, T. Hudson, M. Lamoureux
and B. Solel, \textit{Meet irreducible ideals in direct limit algebras}
Math. Scand., \textbf{87}, 2000, 27--63.

\bibitem{DHu} A. Donsig and T. Hudson,
\textit{The lattice of ideals
of a triangular AF algebra}, J.\ Funct.\ Anal.\ \textbf{138}, 1996, 1--39.

\bibitem{DHK} A. Donsig, T. Hudson and E. Katsoulis,
\textit{Algebraic isomorphisms of limit algebras},
Trans.\ Amer.\ Math.\ Soc.\ \textbf{353}, 2001, 1169--1182

\bibitem{DPP} A. Donsig, D. Pitts, S. Power,
\textit{Algebraic isomorphisms and spectra of triangular limit algebras},
Indiana Univ.\ Math.\ J.\ \textbf{50}, 2001, 1131--1147.

\bibitem{Erd} J. Erdos,
\textit{ Operators of finite rank in nest algebras}
J.\ London Math.\ Soc.\ \textbf{43}, 1968, 391--397.

\bibitem{ErdP} J. Erdos and S. Power,
\textit{ Weakly closed ideals of nest algebras}
J.\ Operator Theory \textbf{7}, 1982, 219--235.

\bibitem{Gl} J. Glimm,
\textit{Type I $\ca$-algebras},
Ann.\ Math.\ \textbf{73}, 1961, 572--612.

\bibitem{HPP00} A. Hopenwasser, J. Peters, S. Power,
\textit{Nest Representations of TAF Algebras},
Canad. J. Math. \textbf{52}, 2000, 1221--1234.

\bibitem{KaK} E. Katsoulis and D. Kribs,
\textit{Isomorphisms of algebras associated with directed graphs},
e-print arxiv math.OA/0309363, preprint, 2003.

\bibitem{KaM} E. Katsoulis and R. Moore,
\textit{On compact operators in certain reflexive operator algebras},
J.\ Operator Theory \textbf{25}, 1991, 177--182.

\bibitem{KrP} D.W. Kribs, S.C. Power,
\textit{Free semigroupoid algebras},
preprint, 2003.

\bibitem{La93} M. Lamoureux,
\textit{Nest Representations and Dynamical Systems},
J. Funct. Anal. \textbf{114}, 1993, 345--376.

\bibitem{La96} M. Lamoureux,
\textit{Ideals in some continuous nonselfadjoint crossed product algebras},
J. Funct. Anal. \textbf{142}, 1996, 211--248.

\bibitem{La97} M. Lamoureux,
\textit{Some Triangular AF Algebras},
J. Operator Theory \textbf{37}, 1997, 91--109.

\bibitem{LaS} D. Larson and B. Solel,
\textit{ Structured triangular limit algebras}
Proc. London Math. Soc. \textbf{75}, 1997, 177--193.

\bibitem{Muh} P.S. Muhly,
\textit{A finite dimensional introduction to operator algebra},
A. Katavolos (ed.), Operator Algebras and
Application, Kluwer Academic Publishers, 1997, 313--354.

\bibitem{MuS} P. Muhly and B. Solel,
\textit{Subalgebras of groupoid
$C\sp *$-algebras},
J.\ Reine Angew.\ Math. \textbf{402}, 1989, 41--75.

\bibitem{OPe}  J. Orr and J. Peters,
\textit{Some representations of TAF algebras},
Pacific.\ J.\ Math. \textbf{167} 1995, 129--161.

\bibitem{Paul} V. I. Paulsen,
\textit{Completely Bounded Maps and Dilations},
Longman Scientific, New York, Wiley, 1986.

\bibitem{Ped} G. K. Pedersen,
\textit{C$^*$-Algebras and their automorphism groups},
London Mathematical Society Monograph, Academic Press, 1979.

\bibitem{PPW} J. Peters, Y.T. Poon and B. Wagner,
\textit{Triangular AF algebras}, J.\ Operator Theory \textbf{23} 1990,
81--114.

\bibitem{Po} S. Power,
\textit{Classification of tensor products of triangular
operator algebras},
Proc.\ London Math.\ Soc.\ \textbf{61}, 1990, 571--614.

\bibitem{Ringr} J. R. Ringrose,
\textit{On some algebras of operators},
Proc.\ London Math.\ Soc.\ \textbf{15}, 1965, 61--83.

\bibitem{SV} \c S. Str\u atil\u a and D.V. Voiculescu,
\textit{Representations of AF-algebras and of the group $U(\infty)$},
Springer Lect.\ Notes Math.\ \textbf{486},
Springer Verlag, Berlin, New York, 1975.

\bibitem{Thel91} M. Thelwall,
\textit{Dilation theory for subalgebras of AF algebras},
J. Operator Theory \textbf{25}, 1991, 275--282.

\end{thebibliography}
\end{document}